\documentstyle[12pt]{article}

\newif\ifgothicavailable

\ifgothicavailable
\newfont{\goth}{eufm10 scaled 1200}   
\else
\newcommand{\goth}{\bf}
\fi

\newcommand{\N}{\mbox{\goth N}}

\advance\textwidth0.5cm
\advance\textheight0.5cm

\newcommand{\fin}[1]{[#1]^{<\omega}}
\newcommand{\fink}[2]{[#2]^{#1}}

\renewcommand{\P}{{\cal P}}

\newcommand{\supp}{\mathop{{\rm supp}}}

\def\cl#1:{\paragraph{#1:}}
\newcommand{\ii}[1]{\iii{f}{#1}}
\newcommand{\iii}[2]{{#1}^{-1}[\{#2\}]}

\newenvironment{ITEMIZE}%
	{\begin{list}{$*$}{\labelwidth1.6cm\leftmargin2cm\rightmargin1cm}}%
	{\end{list}}

\newcommand{\proof}{\par\noindent{\it Proof:} }

\makeatletter
\def\enumerate{\ifnum \@enumdepth >3 \@toodeep\else
      \advance\@enumdepth \@ne
      \edef\@enumctr{enum\romannumeral\the\@enumdepth}\list
      {\csname label\@enumctr\endcsname}{\parsep0pt\itemsep0pt
\usecounter
	{\@enumctr}\def\makelabel##1{\hss\llap{##1}}}\fi}

\def\itemize{\ifnum \@itemdepth >3 \@toodeep\else \advance\@itemdepth \@ne
\edef\@itemitem{labelitem\romannumeral\the\@itemdepth}%
\list{\csname\@itemitem\endcsname}{\parsep0pt\itemsep0pt
\def\makelabel##1{\hss\llap{##1}}}\fi}

\makeatother

\def\xx#1 {%
\newtheorem{#1}[thm]{#1}}
\def\yy#1 {%
\newenvironment{#1}{\begin{roman-#1}\rm}{\end{roman-#1}}
\newtheorem{roman-#1}[thm]{#1}}

\yy Fact
\xx Theorem
\xx Lemma
\yy Definition
\yy Remark
\yy Notation
\yy {Fact and Definition}

\begin{document}

\nocite{Degen}
\nocite{Diel}
\nocite{Fraenkel}
\nocite{HowardYorke}
\nocite{Jech}
\nocite{JechSochor}
\nocite{Kunen}
\nocite{Levy:58}
\nocite{SpisiakVojtas}
\nocite{Tarski:24}
\nocite{Tarski:30}
\nocite{Truss:74}

\author{Martin Goldstern}
\title{A Note on Superamorphous Sets and
\\ 
Dual Dedekind-Infinity}

\date{April 3, 1995}
\maketitle

\section{Introduction}

In the absence of the axiom of choice there  are several 
possible nonequivalent  ways of translating the intuitive idea of
``infinity'' into a mathematical definition.   

In 
\cite{Tarski:24}, Tarski 
 investigated some natural infinity notions notions,
and his research 
was continued by Levy \cite{Levy:58}, Truss \cite{Truss:74}, 
Spi\v siak and Vojtas \cite{SpisiakVojtas}, Howard and Yorke
\cite{HowardYorke} and others. 

The most prominent definitions of finiteness are the
 {\em true finiteness}
 (equipotent to a bounded set of natural numbers), which is
equivalent to 
{\em  Tarski-finiteness} (every family of subsets has a maximal
element), and the much weaker {\em Dedekind-finiteness}
 (not equipotent to
any proper subset).

In a recent survey paper \cite{Degen}, Degen recalled  the notions of
``weak Dedekind infinity'' and ``dual Dedekind infinity'' (see below
for definitions) and asked if
they were in fact equivalent.  A forcing construction of a weakly
Dedekind  set which
is not dually Dedekind is given by Truss in \cite{Truss:74}, and also
announced by Diel \cite{Diel}. 

We give here an alternate and more elementary  
 construction of such a set, which does not use the methods of forcing
but instead relies on (the consistency of) the existence of a
superamorphous set. 

Answering another question posed in \cite{Degen}, we also show that
``inexhaustibility'' is not a notion of infinity unless the axiom of
choice is assumed. 

\bigskip

\begin{Notation}
\begin{enumerate}
\item For any set $A$ we let $A+1$ be any set of the form $A\cup \{a\}$
with $a\notin A$. 
\item $A \le B$ means that there is an injective (i.e., 1-to-1)
function from $A$ into~$B$. 
\item $A\le^* B $ means that there is a surjective function from $B$
onto~$A$. 
\item $\P(A)$ is the power set of~$A$. 
\item $\fink k A$ is the set of $k$-element subsets of $A$, for $k\in
	\omega $.  ($\omega = \{0,1,2,\ldots\}$.)
\item For a finite set $s$ we let $|s|$ be the cardinality of $s$. 
\item $\fin A = \bigcup_{k\in \omega} \fink k A$ is the set of 
	finite subsets of $A$. 
\end{enumerate}
\end{Notation}

\begin{Definition}
We call an infinite set $A$
\begin{itemize}
\item Dedekind infinite (D-infinite) if there is a injective nonsurjective
map from $A$ into~$A$. 
\item dually Dedekind infinite (dD-infinite), if  there is a surjective
noninjective map from $A$ onto~$A$.
\item weakly Dedekind infinite (wD-infinite), if there is a surjective map
from $A$ onto the natural numbers. 
\end{itemize}
\end{Definition}

\begin{Remark}
D-finite sets are the sets that are finite in the IV-th sense  
\cite{Tarski:24}, \cite{Levy:58},  \cite{SpisiakVojtas} or $\Delta$-infinite
in \cite{Truss:74}.   It follows easily from the definitions that 
\begin{ITEMIZE}
\item \quad  $A$ is D-infinite \quad   iff  \quad $ A+1 \le A$ \quad iff
\quad  $ \omega \le A$ 
\end{ITEMIZE}

dD-finite sets are those in  $\Delta_5$ in \cite{Truss:74}, or ``not
strongly Dedekind finite'' in the sense of \cite{Diel}.  Trivially, 
\begin{ITEMIZE}
\item \quad $A$ is dD-infinite \quad iff \quad $  A+1\le^* A$ 
\end{ITEMIZE}

wD-finite sets are those sets that are finite in the III-rd sense in
\cite{Tarski:24} etc, $\Delta_4$ in \cite{Truss:74}, and ``almost
finite'' in \cite{Diel}.    It is
well-known and not hard to show that 
\begin{ITEMIZE}
\item \quad 
$A$ is wD-infinite \quad iff \quad  $   \omega \le^* A$ \quad iff \quad
$\omega \le \P(A)$
\end{ITEMIZE}

We have D-infinite $\Rightarrow $ dD-infinite $\Rightarrow$ wD-infinite, and none of
these implications can be reversed in ZF.  For example, if $A$ is
amorphous (see \ref{amorphdef})
then the set of injective finite sequences from $A$ is a set
that is dD-infinite but not D-infinite. 
\end{Remark}

The following definition appears to be new:

\begin{Definition}
A set $A$ is {\em wD*-infinite} if there is a finite-to-one map from a
subset of $A$ onto $\omega$ (or equivalently, onto an infinite subset
of $\omega$). 
\end{Definition}

\begin{Fact}\label{wdstar}
For any set $A$ the following are equivalent: 
\begin{enumerate}
\item $A$ is wD*-infinite.
\item $\fin A$ is D-infinite.
\item $\fin A$ is wD*-infinite.
\item There is a sequence $(A_i:i\in \omega)$, $A_i\in \fin A$, 
	$\bigcup_{i\in \omega} A_i$ infinite. 
\item There is a strictly increasing 
	sequence $(A_i:i\in \omega)$, $A_i\in \fin A$, 
	$\bigcup_{i\in \omega} A_i$ infinite. 
\end{enumerate}
\end{Fact}
\proof  The implications 
(4)
$\Rightarrow$
(5)
$\Rightarrow$
(1) 
$\Rightarrow$
(2)
$\Rightarrow$
(3)
are trivial.   For (3) $\Rightarrow$ (4), let $f$ be a map with domain
$\subseteq \fin A$ and range $= \omega$ such that the preimage of any
natural number is finite.   For $n\in \omega $ let 
$B_n:=\{s : f(s) = n\}$, and let $A_n:= \bigcup_{s\in B_n} s$. 
Then the domain of $f$ is contained in $\fin {\bigcup_n A_n}$, so
$\bigcup_n A_n$ must be infinite, whereas all the sets $A_n$ are
finite.

\begin{Remark}
D-infinite $\Rightarrow$ wD*-infinite $\Rightarrow$ wD-infinite, but
none of 
these implications can be reversed.  Moreover, there is (in ZF) no
implication between dD and wD*: The set of finite injective sequences
from an amorphous set is dD-infinite but not wD*-infinite, and the set
$U$ in the model $\N_2$ of \cite{Truss:74} (a union of countably many
``pairs of socks'') is clearly wD*-infinite but not
dD-infinite. 

Slightly more generally, if $(U_n:n\in \omega)$ is a sequence of pairwise
disjoint unordered pairs such that $\bigcup_n U_n$ is Dedekind-finite,
(or equivalently, 
such that no infinite subsequence
$(U_{n_i}:i\in \omega)$ has a choice function) 
then it is easy to see that 
$\bigcup_n U_n$ is  dD-finite. 
\end{Remark}

\section{Superamorphous sets}

\begin{Definition}\label{amorphdef}
A set $A$ is called amorphous, if $A$ is infinite and 
all its subsets are either finite or cofinite. 
\\
A set 
$A$ is called superamorphous, if $A$ is infinite and for all $k$, 
all subsets of $A^k$ are first order definable from finitely many
parameters in the language of equality. 
\end{Definition}

\begin{Remark}  It is consistent with ZF 
that there are superamorphous sets.
\\
 The consistency of the existence of superamorphous
sets with a nonwellfounded set theory or with ZFU (=ZF with urelements)
can already be seen in the basic Fraenkel-Mostowski model (see
\cite{Fraenkel} and 
\cite[Exercise IV.24]{Kunen}).  The  consistency with ZF then follows
from the Jech-Sochor theorem (see \cite[Theorem 47]{Jech},
\cite{JechSochor}). 
\end{Remark}

Using a superamorphous set, 
we give an elementary example of a wD-infinite  set that is not
dD-infinite by the following theorem:

\begin{Theorem} If $A $
is a superamorphous set, then the set $\P(A)$ is wD-infinite but not
dD-infinite. 
\end{Theorem}

Before we prove the theorem, we collect some facts about~$\fin A$. 
The following fact will allow us to work with $\fin A$ instead of
$\P(A)$. 

\begin{Fact}\label{A+1}
Let $A$ be (super)amorphous. Then also 
$A+1$ is (super)amorphous, and there is a bijection from
$\P(A)$ onto $\fin {A+1}$. 
\end{Fact}

We will need a ``normal form'' for subsets of $\fink k A$. 
\begin{Definition}
Fix a set $A$ and a natural number~$k$.  Let 
$p,q \in \fin A$. We let 
$$ A^{+p-q}(k):= \{s\in \fink k A: p \subseteq s, q\cap
	s=\emptyset\}$$
\end{Definition}

\begin{Fact}\label{nonempty}
Assume that $A$ is infinite.  Then:
\begin{enumerate}
\item $A^{+p-q}(k)\not=\emptyset$ iff $p\cap q = \emptyset$ 
			and $|p|\le k$. 
\item $A^{+p-q}(k) \cap A^{+p'-q'}(k)= 
			A^{+p\cup p'\, -\, q\cup q'}(k)$. 
\item $ A^{+p-q}(k)$ is infinite iff  $p\cap q = \emptyset$ 
			and $|p|<  k$. 
\item If $A^{+r-s}(k)$ is infinite, then there are $p,p'\in 
			A^{+r-s}(k)$ with $|p\cup p'| = k+1$. 
\end{enumerate}
\end{Fact}

\begin{Fact and Definition}
Let $A$ be an infinite set. 
Assume that $B \subseteq A^k$ is definable.  Then the set 
$$ B^*:= \{ s\in \fin A : \hbox{$B$ is definable from parameters in $s$}
\}$$
has a smallest element $\bigcap B^*$.  We call this smallest element the
``support'' of $B$, $\supp(B)$.
\end{Fact and Definition}

\begin{Fact}[normal form]\label{normalform}
Let $C \subseteq \fink k A $ and assume that $\bar C:=
\{(x_1, \ldots, x_k): \{x_1, \ldots, x_k\} \in C\}$ 
is definable.  Let 
$$ H_C:= \{(p,q): p \cup q = \supp(\bar C),\, p \cap q=\emptyset, \,
	A^{+p-q}(k) \subseteq C \}$$
Then $C = \bigcup_{(p,q)\in  H_C} A^{+p-q}(k)$. 
(Notice that $H_C$ is a finite set.) 
\end{Fact}
\proof  If $s\in C$, then let $p:= s\cap \supp(\bar C)$, 
$q:= \supp(\bar C) \setminus p$.  Then 
$s \in A^{+p-q}(k) \subseteq C$.


 \begin{Lemma}\label{nobigunion}
Let $A$ be amorphous.  Then: 
 \begin{enumerate}
 \item \label{asubi}
	$\fin A$ is wD*-infinite. 
 \item \label{asupi}
	If  $B \subseteq \fin A$ is  infinite, then 
	there is some $k$ such that $B \cap \fink k A
	$ is infinite.
 \end{enumerate}
 \end{Lemma}
 \proof (1) follows from \ref{wdstar}.  If (2) were false, then the
map $s \mapsto |s|$,  would witness that $B$ and hence also $\fin A $
is wD*-infinite, which is impossible.

\begin{Definition}
For $R \subseteq X \times Y$, $x\in X$, $y\in Y$ we let 
$R^x:= \{z\in Y: (x,z)\in R\}$, and 
similarly we let $R_y:= \{z\in X: (z,y)\in R\}$.
\end{Definition}

\begin{Lemma}\label{rkl}
Assume that $A$ is superamorphous, $k,l \in \omega$, $R \subseteq 
\fink k A  \times \fink l A$ is infinite, and for all $(p,q)\in R$ we
have $p\cap q = \emptyset$. \\
Then either there is some $p\in  \fink k A$ such that $R^p$ is
infinite, or there is some $q\in  \fink l A$ such that $R_q$ is
infinite (or both). 
\end{Lemma}

\subsubsection*{Proof of the theorem}

By fact \ref{A+1} it is enough to show that for any superamorphous set
$A$, every 
surjective map from $\fin A $ onto $ \fin A$ is injective.

Let $f: \fin A \to \fin A$ be onto. We will show 
that for all $s \in \fin A$ there is $\exists n>0 $ such that
$f^n(s)=s$.          
(This easily yields that $f$ is injective.) 

We will proceed indirectly and assume that
\begin{ITEMIZE}
\item[$(\ast)$] there is $s$ such that 
 $\forall n\in \omega\,\, f^n(s)\not=s$. 
\end{ITEMIZE}

\cl Claim 1: If $\forall n\in \omega\,\, f^n(s)\not=s$, then there is
$s'$ with $\ii{s'} $ infinite, and $\exists n\, f^n(s')=s$. 
\proof Let $X:= \fin A $. If there is no $s'$ as claimed, then for
each $i\in \omega $ the set 
$$X_i:= \{t\in X: f^i(t) = s\}$$
would be in $\fin X$, and all these sets would be (nonempty and) distinct.
Hence the function $i\mapsto X_i$ 
would witness that $\fin X$ is
D-infinite.  Applying \ref{wdstar} twice we would get that $X$ and
hence $A$  is
wD*-infinite, a contradiction.

\cl Claim 2: If for all $n$ we have $f^n(s)\not=s$, then there are
infinitely many $s'$ as in claim~1. 
\proof Apply Claim 1 repeatedly.

\cl Claim 3: If $\ii{s'}$ is infinite, then there is some $k$ 
such that $\ii{s'} \cap \fink k A$ is infinite. 
\proof By \ref{nobigunion}. 

\cl Claim 4: For some $k$, there are infinitely many $s\in \fin A$
such that the set $\ii s \cap \fink k A$ is infinite. 
\proof 
Let $B_k:= \{s': \ii{s'} \cap \fink k A \mbox{ is infinite}\}$.  By
claims 2~and~3, $\bigcup_{k\in \omega } B_k$ is infinite, so 
 some $B_k$ must be  infinite.

\cl Claim 5: For some $k$ there is an infinite set $R \subseteq
\fin A \times \fin A $ such that for all $(p,q)\in R$ the set 
$A^{+p-q}(k)$ is infinite, but any two such sets are disjoint. 
\proof Fix $k$ as in  claim 4. Replace each infinite set of the form 
$\ii s \cap \fink k A$ by those components $A^{+p-q}(k)$ in its normal
form (see \ref{normalform}) that are infinite. 

\cl Claim 6: There are $n_1$ and $n_2$ such that we can find $R$ as in
claim 5, additionally satisfying 
$|p | = n_1$ and 
$|q | = n_2$ for all $(p,q) \in R$.  (By \ref{nonempty} we must then
have $n_1 < k$.)
\proof  Easy.

Now let $R$ be as above. 
$R$ satisfies the assumptions of lemma \ref{rkl}.   We distinguish two
cases:  
\cl Case 1: For some $p\in \fink {n_1} A$, $R^p$ is infinite. 
In particular this means that there are $q\not= q'$ in~$R^p$.
But then $\emptyset \not= 
A^{+p-(q\cup q')}(k) \subseteq A^{+p-q}(k) \cap
A^{+p-q'}(k)$, a contradiction. 

\cl Case 2: 
 For some $q\in \fink {n_1} A$, $R_q$ is infinite. 
By fact \ref{nonempty} we can find $p, p'\in R_q$ with 
$|p\cup p'| =n_1+1 \le k$.  So  we have 
$\emptyset \not= A^{+p\cup p'-q}(k) \subseteq A^{+p-q}(k)\cap 
A^{+p'-q}(k)$, again a  contradiction to what
we have found in claim 5.

\section{Inexhaustibility}
Finally we answer another question from Degen's paper: 
\begin{Definition}
An  set $A$ is called ``inexhaustible'' if $A$ contains more than one
element, and for any
decomposition $A= B \cup C$ we have that $A$ can be injectively
mapped into $B$ or~$C$.  
\end{Definition}

It is clear that if AC holds, then inexhaustible sets are exactly the
infinite sets.  However, without AC 
 it is not clear if
supersets of inexhaustible sets are also inexhaustible.    We show
that this property actually characterizes AC: 

\begin{Theorem}
Assume that every set that contains an inexhaustible set is itself
inexhaustible.  Then the axiom of choice holds. 
\end{Theorem}

If we call a property a ``notion of infinity'' (as in \cite{Degen})
 iff it is closed under
equivalences and supersets, and holds for $\omega$ but no finite set,
then we can rephrase the above theorem as follows:

\begin{quote}
``Inexhaustibility'' is a notion of infinity only if the axiom of
choice holds, or equivalently, only if it coincides with true 
infinity. 
\end{quote}

\proof
 We will show that every set can be well-ordered.   This is
clear for finite sets, so consider some infinite~$B$. 

Let $\kappa$ be the least ordinal number such that $B$
cannot be mapped onto $\kappa$.
Let $A$ be the disjoint union of $\kappa$ and $B$: 
$$ A = B +  \kappa $$
Clearly $\omega 
\subseteq \kappa$, hence $A$ contains an inexhaustible set and is
therefore itself inexhaustible. \\
 So either we have $B + \kappa \le B$ or $B+ \kappa  \le \kappa$.  The
first alternative 
would imply $\kappa \le B$, hence $\kappa \le^* B$ which is impossible
by the definition of~$\kappa$. 
\\
Hence we have $B + \kappa  \le \kappa$, so also $B \le \kappa$.
Hence $B$ can be well-ordered.

\begin{Remark}  The use of ordinals, and hence of the replacement axiom, can
be avoided.  Instead of $\kappa$, take the set of all quasiorders on
$B$ such that the quotient p.o.\ is a well-ordering.  This set is
naturally well-ordered and can play the role of $\kappa$ in the above
proof.
\end{Remark}

\bibliographystyle{plain}

\bigskip
\vfill

\hrule 
\bigskip

\bigskip

\noindent
Martin Goldstern, 
Department of Algebra, 
Technische Universit\"at, \\
Wiedner Haupstra{\ss}e 8-10/118.2, 
A-1040 Wien, Austria, Europe
\\
{\tt martin.goldstern@tuwien.ac.at}
\eject

\end{document}